\documentclass{elsart}
\usepackage{amssymb}
\usepackage{amsmath}
\usepackage{lscape}
\journal{Journal of Algebra}

\newcommand\A{{\mathcal A}}
\newcommand\B{{\mathcal B}}
\newcommand\U{{\mathcal U}}
\newcommand\N{{\mathcal N}}
\newcommand\vng{{\mathcal N}G}
\newcommand\ug{{\mathcal U(G)}}
\newcommand\ann{\mathrm{ann}}
\newcommand\Qrmax{Q^r_{\mathrm{max}}}
\newcommand\Qmax{Q_{\mathrm{max}}}

\newcommand\lge{l^2(G)}

\begin{document}
\begin{frontmatter}

\title{$\ast$-Clean Rings; Some Clean and Almost Clean Baer $\ast$-rings and von Neumann Algebras}

\author{Lia Va\v s}
\address{Department of Mathematics, Physics and Statistics\\ University of the Sciences
in Philadelphia\\ Philadelphia, PA 19104, USA}

\ead{l.vas@usp.edu}

\begin{abstract}
A ring is clean (resp. almost clean) if each of its elements is the sum of a unit (resp. regular element) and an
idempotent. In this paper we define the analogous notion for $\ast$-rings: a
$\ast$-ring is $\ast$-clean (resp. almost $\ast$-clean) if its every element is the sum of a unit (resp. regular element) and a
projection. Although $\ast$-clean is a stronger notion than clean, for some $\ast$-rings we demonstrate that it is more natural to use.

The theorem on cleanness of unit-regular rings from \cite{Camillo}
is modified for $\ast$-cleanness of $\ast$-regular rings that are abelian (or reduced or Armendariz). Using this result, it is shown that all finite, type $I$ Baer $\ast$-rings that satisfy certain axioms (considered in \cite{Be} and \cite{Lia2}) are almost $\ast$-clean. In particular, we obtain that all finite type $I$ $AW^{\ast}$-algebras (thus all finite type $I$ von Neumann algebras as well) are almost $\ast$-clean.  We also prove that for a Baer $\ast$-ring satisfying the same axioms, the following properties are equivalent: regular, unit-regular, left (right) morphic and left (right) quasi-morphic. If such a ring is finite and type $I$, it is $\ast$-clean. Finally, we present some examples related to group von Neumann algebras and list some open problems.
\end{abstract}


\end{frontmatter}

\section{Introduction}

The motivation for this paper came from a question posed by T.Y. Lam at the Conference on Algebra and Its Applications held at Ohio University, Athens, OH in March, 2005. Lam asked which von Neumann algebras are clean as rings.

A ring is {\em clean} if its every element can be written as the sum
of a unit and an idempotent. Clean rings were defined by W. K. Nicholson in late 1970s in relation to exchange rings (see \cite{Nich}) and have been attracting attention ever since. Clean rings are, in a way, an additive analogue of unit-regular rings: in a unit-regular ring,  every element can be written as the product of a unit and an idempotent, for clean rings ``the product'' in this condition changes to ``the sum''. Thus, it was not a coincidence that the question of the exact relationship between the classes of unit-regular and clean rings was of particular interest. In \cite{Camillo}, the authors characterize unit-regular rings as clean rings in which every element $a$ can be decomposed as the sum of a unit $u$ and an idempotent $e$ with $aR\cap eR=0$ (Camillo-Khurana Theorem). Other examples of clean rings include local rings, semiperfect rings and right artinian rings. A nice overview of commutative clean rings is given in the first section of \cite{Warren_neat}.

A ring is {\em almost clean} if its every element can be written as a sum of a regular element (neither a left nor a right zero-divisor) and an idempotent. Up to date, almost cleanness has been considered mostly for commutative rings. The concept was introduced in \cite{Warren_CX}. In \cite{Warren_CX}, it is shown that a commutative Rickart ring is almost clean.  In this paper, we shall exhibit various classes of almost clean rings that are not necessarily commutative.

Going back to Lam's question, let us turn to von Neumann algebras. Since a von Neumann algebra is a $\ast$-ring (i.e. a ring with involution) in some cases it may be more natural to work with projections, self-adjoint idempotents, rather than with idempotents. Many other examples of different matrix rings and various rings of operators naturally come equipped with an involution. For such rings the properties of being $\ast$-regular, Baer $\ast$-ring or Rickart $\ast$-ring take over the roles of regular, Baer or (right or left) Rickart respectively (see, for example, \cite{Be}). Also, as S. K. Berberian points out in \cite{Be}, ``the projections are vastly easier to work with than idempotents'' so one should utilize this opportunity when dealing with a $\ast$-ring. Because of this, we define the following more natural concept of cleanness for $\ast$-rings.

\begin{defn}
A $\ast$-ring is {\em $\ast$-clean (resp. almost $\ast$-clean)} if its every element can be written as the
sum of a unit (resp. regular element) and a projection.
\end{defn}
Clearly, if a $\ast$-ring is (almost) $\ast$-clean it is also (almost) clean.

This definition is further justified by the following observation: if a unital ring $R$ embeds in a clean (resp. $\ast$-clean) ring $S$ with the same idempotents (resp. projections), then $R$ is almost clean (resp. almost $\ast$-clean). Every finite von Neumann
algebra $\A$ embeds in a unit-regular ring $\U$ (the algebra of affiliated operators). Moreover the ring $\U$ has the same projections as $\A.$ Since unit-regular rings are clean by the Camillo-Khurana Theorem (\cite[Theorem 1]{Camillo}), $\U$ is clean. Thus, we may hope that the cleanness of $\U$ will be preserved in $\A$ at least to the extent of almost cleanness. The problem is that the Camillo-Khurana Theorem does not guarantee that one can decompose an element in a unit-regular $\ast$-ring as the sum of a unit and a projection, just as the sum of a unit and an idempotent. We would need the stronger version to pull this decomposition from $\U$ back down to $\A.$ However, we shall demonstrate that for a class of von Neumann algebras the $\ast$-clean decomposition in $\U$ is possible and that it yields almost $\ast$-cleanness of $\A.$

The paper is organized as follows. In Section \ref{section_preliminaries}, we recall some known concepts and results. In Section \ref{section_abelian_*-clean}, we prove some preliminary results and basic properties of $\ast$-clean rings. We also modify the Camillo-Khurana Theorem to fit the $\ast$-ring setting and prove that abelian (or reduced or Armendariz) $\ast$-regular rings are $\ast$-clean (Proposition \ref{clean_regular}). We also show that every $\ast$-ring that embeds in a $\ast$-clean regular $\ast$-ring with same projections is almost $\ast$-clean (Proposition \ref{embedding_proposition}).

In Section \ref{section_Baer}, we consider a class of Baer $\ast$-rings for which the involution extends to their maximal right ring of quotients ($\ast$-extendible rings). We recall Baer $\ast$-ring axioms (A1)--(A7) considered in \cite{Be} and \cite{Lia2} and prove a result stating that every type $I_n$ Baer $\ast$-ring satisfying (A2) is almost $\ast$-clean (Theorem \ref{type_In}). Using that, we show that every type $I$ Baer $\ast$-ring that satisfies (A1)--(A6) is almost $\ast$-clean (Theorem \ref{type_If}). In particular, this theorem gives us that all finite, type $I$ $AW^{\ast}$-algebras (thus all finite type $I$ von Neumann algebras as well) are almost $\ast$-clean (Corollary \ref{VNA_clean}). If regularity is also assumed in Theorems \ref{type_In} and \ref{type_If} and Corollary \ref{VNA_clean}, we obtain $\ast$-cleanness. In Section \ref{section_Baer}, we also prove that for a Baer $\ast$-ring that satisfies (A1)--(A6), the following properties are equivalent: regular, unit-regular, left (right) morphic and left (right) quasi-morphic (Corollary \ref{morphic}).

In Section \ref{section_VNAs}, we turn to $AW^{\ast}$ and von Neumann algebras. Since every finite $AW^{\ast}$-algebra of type $I$ is almost $\ast$-clean (as Corollary \ref{VNA_clean} shows), we wonder if there are examples of almost $\ast$-clean $AW^{\ast}$-algebras that are not type $I.$ We show that a $C^{\ast}$-sum of finite type $I$ $AW^{\ast}$-algebras is almost $\ast$-clean (Proposition \ref{AW_C*-sums}). Moreover, a $C^{\ast}$-sum  of finite type $I$ regular $AW^{\ast}$-algebras is $\ast$-clean. Using this fact, we consider a group von Neumann algebra that is not of type $I$ but is $\ast$-clean. We also consider other group von Neumann algebras (Example \ref{Group_Example}) that provide examples of some other possible situations. We conclude the paper with a list of open problems in Section \ref{section_questions}.

\section{Preliminaries}
\label{section_preliminaries}

An associative unital ring $R$ is a {\em $\ast$-ring} (or {\em ring with
involution}) if there is an operation $\ast: R\rightarrow R$ such
that
\[(x+y)^{\ast}=x^{\ast}+y^{\ast},\;\;\;(xy)^{\ast}=y^{\ast}x^{\ast},\;\;\;(x^{\ast})^{\ast}=x\;\;\;\mbox{ for all }x,y\in R.\]
If a $\ast$-ring $R$ is also an algebra over $k$ with involution $\ast,$ then $R$
is a {\em $\ast$-algebra} if $(ax)^{\ast}=a^{\ast}x^{\ast}$ for $a\in k,$ $x\in R.$
An element $p$ of a $\ast$-ring $R$ is called a {\em projection} if $p$
is a self-adjoint ($p^{\ast}=p$) idempotent ($p^2=p$). There is an equivalence relation on the set of projections of a $\ast$-ring $R$ defined by $p\sim q$ iff $x^{\ast}x=p$ and $xx^{\ast}=q$ for some
$x\in R.$

Recall that a ring is {\em right Rickart} (also called right p.p.-ring) if the right annihilator of each element is generated by an idempotent. A ring is said to be {\em Rickart $\ast$-ring} if the right annihilator of each element is generated by a projection. Moreover, it can be shown that such  projection is unique. Since $\ann_l(x)=(\ann_r(x^{\ast}))^{\ast}$ the property of being Rickart $\ast$-ring is left/right symmetric. In every Rickart $\ast$-ring, the involution is proper ($x^{\ast}x=0$ implies $x=0$, see \cite[Proposition 2, p. 13]{Be}). The projections in a Rickart $\ast$-ring form a lattice (\cite[Proposition 7, p. 14]{Be}) with $p\leq q$ iff $p=pq$ (equivalently $p=qp$).

Every element $x$ of a Rickart $\ast$-ring $R$ determines a unique
projection $p$ such that $xp=x$ and $\ann_r(x)=\ann_r(p)=(1-p)R$
and a unique projection $q$ such that $qx=x$ and
$\ann_l(x)=\ann_l(q)=R(1-q).$ In this case, $p$ is called the {\em right projection of $x$} and is denoted by RP$(x),$ and $q$ is the {\em left projection of $x$} and is denoted by LP$(x).$

A $\ast$-ring $R$ is $\ast$-regular if any of the following equivalent
conditions hold (see \cite[Proposition 3, p. 229]{Be}): (1) $R$ is regular and Rickart $\ast$-ring; (2) $R$ is regular and the involution is proper; (3) $R$ is regular and for every $x$ in $R$ there exists a
projection $p$ such that $xR=pR.$

A {\em Rickart $C^{\ast}$-algebra} is a $C^{\ast}$-algebra (complete normed
complex algebra with involution such that $||a^{\ast}a||=||a||^2$) that
is also a Rickart $\ast$-ring.

A ring is {\em Baer} if the right (equivalently left) annihilator of every nonempty subset is generated by an idempotent. A {\em Baer $\ast$-ring} is a $\ast$-ring $R$ such that the right annihilator of every nonempty subset is generated by a projection. A $\ast$-ring is a Baer $\ast$-ring if and only if it is a Rickart $\ast$-ring and the lattice of projections is complete (\cite[Proposition 1, p. 20]{Be}).

An {\em $AW^{\ast}$-algebra} is a $C^{\ast}$-algebra that is a Baer $\ast$-ring. If $H$ is a Hilbert space and $\B(H)$ the algebra of bounded operators on $H,$ then $\B(H)$ is an $AW^{\ast}$-algebra. If $\A$ is a unital $\ast$-subalgebra of $\B(H)$ such that $\A= \A''$ where $\A'$ is the commutant of $\A$ in $\B(H),$ then $\A$ is called a {\em von Neumann algebra.} Equivalently, $\A$ is a von Neumann algebra if it is a unital $\ast$-subalgebra of $\B(H)$ that is closed with respect to weak operator topology. A von Neumann algebra is an $AW^{\ast}$-algebra (\cite[Proposition 9]{Be}) while the converse is not true.

\section{$\ast$-clean and almost $\ast$-clean rings}
\label{section_abelian_*-clean}

In this section, we turn to $\ast$-rings. Recall that a ring is abelian if every idempotent is central. The $\ast$-version of this concept is as follows.

\begin{defn}
A $\ast$-ring is {\em $\ast$-abelian} if every projection is central.
\end{defn}

Clearly, a $\ast$-ring that is abelian is also $\ast$-abelian. The following lemma proves that for a Rickart $\ast$-ring, the conditions of abelian and $\ast$-abelian are equivalent and gives some further characterization of abelian right Rickart rings. Also, recall that a ring is said to be {\em strongly clean} if its every element can be written as the sum of a unit and an idempotent that commute. Clearly, a strongly clean ring is clean. If a ring is abelian, then clean implies strongly clean as well.
We shall say that a ring is {\em almost strongly clean} if its every element can be written as the sum of a regular element and an idempotent that commute.

Let us define a $\ast$-ring to be {\em (almost) strongly $\ast$-clean} if its every
element can be written as the sum of a unit (resp. regular element) and a projection that
commute.

\begin{lem}
\begin{enumerate}
\item If $R$ is a $\ast$-abelian Rickart $\ast$-ring, then every idempotent is a projection.

\item If $R$ is a Rickart $\ast$-ring, then $R$ is abelian iff $R$ is $\ast$-abelian.

\item Let $R$ be a ($\ast$-)abelian Rickart $\ast$-ring. Then the following conditions are equivalent:

\begin{itemize}
\item[(i)] The ring $R$ is $\ast$-clean.
\item[(ii)] The ring $R$ is  clean.
\item[(iii)] The ring $R$ is strongly ($\ast$)-clean.
\end{itemize}

\item If $R$ is an abelian Rickart $\ast$-ring, then $R$ is almost clean, almost $\ast$-clean and almost strongly ($\ast$)-clean.
\end{enumerate}
\label{abelian_lemma}
\end{lem}
\begin{pf}
(1) Let $e$ be an idempotent in a  $\ast$-abelian Rickart $\ast$-ring $R$. Since $R$ is Rickart $\ast$-ring, there is a projection $p$ such that $eR=\ann_r(1-e)=pR.$ Thus $pe=e$ and $ep=p.$ But since $p$ is central, $e=pe=ep=p.$

Condition (2) follows from (1).

(3) Condition (i) always implies (ii). By part (1), the converse holds for Rickart $\ast$-rings that are abelian (equivalently $\ast$-abelian by (2)). Condition (iii) always implies (ii). The converse clearly holds if all idempotents are central.

(4) By \cite[Proposition 16]{Warren_CX}, a commutative Rickart ring is almost clean. The proof uses \cite[Lemma 2 and Lemma 3]{Endo}. These two lemmas ensure that an element in an abelian right Rickart ring is a product of a regular element and an idempotent. The proof of Proposition 16 uses just that the idempotents are central, not that the ring has to be commutative. Thus, we have that an abelian right Rickart ring is almost clean.

So, if $R$ is an abelian Rickart $\ast$-ring, it is almost clean. But since idempotents are projections by (1), $R$ is also almost $\ast$-clean. Since all the idempotents (projections) are central, $R$ is almost strongly $\ast$-clean as well.
\end{pf}

Before we turn to (almost) $\ast$-cleanness of some specific classes of rings, let us consider the  following proposition proving some basic properties of $\ast$-clean rings. If $R$ is a $\ast$-ring, the ring $M_n(R)$ of $n\times n$ matrices over $R$ has natural involution inherited from $R$: if  $A=(a_{ij}),$ $A^{\ast}$ is the transpose of $(a_{ij}^{\ast}).$ So, $M_n(R)$ is also a $\ast$-ring.

\begin{prop}
Let $R$ be a $\ast$-ring.
\begin{enumerate}
\item If $p$ is a projection such that $pRp$ and $(1-p)R(1-p)$ are $\ast$-clean, then $R$ is $\ast$-clean.

\item If $p_1, p_2,\ldots,p_n$ are orthogonal projections with $1=p_1+p_2+\ldots+p_n,$ and $p_iRp_i$ $\ast$-clean for each $i$, then $R$ is $\ast$-clean.

\item If $R$ is $\ast$-clean, then $M_n(R)$ is $\ast$-clean.
\end{enumerate}
\label{*-clean_matrices}
\end{prop}
\begin{pf}
The proof of (1) follows the proof of \cite[Lemma on p. 2590]{Han_Nicholson} with idempotent in Pierce decomposition of $R$ changed to
projection and every instance of idempotent further in the proof changed to projection. The condition (2) follows from (1) by induction, and (3) follows directly from (2).
\end{pf}

Let us now recall the statement of the Camillo-Khurana Theorem (\cite[Theorem 1]{Camillo}).

\begin{thm}[Camillo-Khurana] A ring $R$ is unit-regular if and only if every
element $a$ of $R$ can be written as $a=u+e$ where $u$ is a unit
and $e$ is an idempotent such that $aR\cap eR=0.$
\label{clean_unit_regular}
\end{thm}

This theorem has the following proposition as a corollary.

\begin{prop} Let $R$ be a $\ast$-ring. If $R$ is $\ast$-regular and abelian, then
every element $a$ of $R$ can be written as $a=u+p$ where $u$ is a
unit and $p$ is a projection such that $aR\cap pR=0.$
\label{clean_regular}
\end{prop}

\begin{pf} The ring $R$ is unit-regular since it is regular and abelian (\cite[Corollary 4.2, p. 38]{Go}). Thus, the condition  equivalent to unit-regularity in Theorem \ref{clean_unit_regular} holds. Since every $\ast$-regular ring is Rickart $\ast$-ring and in every abelian Rickart $\ast$-ring an idempotent is a projection by Lemma \ref{abelian_lemma}, the proposition follows directly from Theorem  \ref{clean_unit_regular}.
\end{pf}

\begin{rem}

{\em
\begin{itemize}
\item[1.] The assumption that $R$ is abelian (let us denote this condition by (1)) in Proposition \ref{clean_regular} can be replaced by any of the conditions below.

\begin{itemize}
\item[(2)] The ring $R$ is {\em semicommutative} if $ab=0$ implies $aRb=0$ for all $a,b\in R$.

\item[(3)] The ring $R$ is {\em symmetric} if $rab=0$ implies $rba=0$ for all $r,a,b\in R$.

\item[(4)] The ring $R$ is {\em Armendariz} if $p(x)q(x)=0$ then $p_iq_j=0$ for all $p(x)=\sum_{i=0}^n p_ix^i$ and $q(x)=\sum_{j=0}^m q_j x^j$ in $R[x]$.

\item[(5)] The ring $R$ is {\em Armendariz of power series type} if $p(x)q(x)=0$ then $p_iq_j=0$ for all $p(x)=\sum_{i=0}^{\infty} p_ix^i$ and $q(x)=\sum_{j=0}^{\infty} q_j x^j$ in $R[[x]]$.

\item[(6)] The ring $R$ is {\em reduced} if it has no nonzero nilpotent elements.
\end{itemize}

In general, $(6)\Rightarrow(3)\Rightarrow (2)\Rightarrow (1),$ $(6)\Rightarrow(5)\Rightarrow (2),$ and $(5)\Rightarrow (4)\Rightarrow (1),$ and each of the implications is strict (see \cite{HKK}, \cite{Kim_Lee_Lee} \cite{Kim_Lee} and \cite{Huh_Lee_Smok}).

There are analogues of conditions (1)--(6) for $R$-modules (see \cite{Rege_Buhpang} and \cite{Agayev_et_al}). In particular, a ring $R$ satisfies one of the properties above if and only if the right module $R_R$ satisfies the corresponding module theoretic
property. \cite[Theorem 2.14]{Agayev_et_al} proves that the module version of the six conditions are all equivalent for a Rickart module. Thus, conditions (1) -- (6) are equivalent for a right Rickart ring.

Since regular rings are Rickart, the assumption that $R$ is abelian in Proposition \ref{clean_regular} can be replaced by any of the conditions (2)--(6). Also, in part (3) of Proposition \ref{embedding_proposition}, the condition that the ring is abelian could be replaced by any of the conditions (2)--(6).

\item[2.] The projection $p$ and unit $u$ from the statement of Proposition \ref{clean_regular} commute since $R$ is abelian so that $R$ is strongly clean as well.

\item[3.] If $R$ is a $\ast$-ring such that every element $a$ of $R$ can be written as $a=u+p$ where $u$ is a
unit and $p$ is a projection such that $aR\cap pR=0,$ then $R$ is unit-regular by Theorem \ref{clean_unit_regular}.
\end{itemize} }
\end{rem}

We will use the following result in the next section.

\begin{prop} Let $R$ be a $\ast$-ring that can be
embedded in a $\ast$-clean regular $\ast$-ring with the same projections as
$R$.
\begin{enumerate}
\item The ring $R$ is almost $\ast$-clean.

\item If $R$ is regular, then $R$ is $\ast$-clean.

\item If $R$ is abelian, then $R$ is almost strongly $\ast$-clean.
\end{enumerate}
\label{embedding_proposition}
\end{prop}

\begin{pf} (1) Let $Q$ be a $\ast$-clean ring in which $R$ is embedded and $a\in R.$ Then $a=u+p$ for a unit $u$ in $Q$ and a projection $p$ in $Q.$ But $p$ is in $R$ by assumption. Thus, $u=a-p$ is in $R$ as well. Since $u$ is a unit in $Q,$ $0=\ann^Q_r(u)\supseteq\ann_r^R(u)$ and the same holds for the left annihilators. Thus $u$ is regular.

(2) If $R$ is regular, then every element is either invertible or a zero-divisor. Thus, the element $u$ from the proof of (1) is invertible and so $R$ is $\ast$-clean.

(3) An abelian almost $\ast$-clean ring is almost strongly $\ast$-clean by Lemma \ref{abelian_lemma}.
\end{pf}

\begin{rem} {\em
Let us note that parts (1) and (3) of Proposition \ref{embedding_proposition} are not valid if ``almost'' is deleted. Note that $\Zset$ (with trivial involution) can be embedded in $\Qset$ with same projections (0 and 1) but $\Zset$ is not (strongly) clean.  }
\end{rem}

\section{$\ast$-clean and almost $\ast$-clean Baer $\ast$-rings}
\label{section_Baer}

In this section, we introduce a class of $\ast$-rings that can be embedded in $\ast$-regular rings with the same projections.  Using Proposition \ref{embedding_proposition}, we shall prove (almost) $\ast$-cleanness of a class of Baer $\ast$-rings.

A $\ast$-ring $R$ is said to be {\em $\ast$-extendible} if its involution can be extended to an involution of its maximal right ring of quotients $\Qrmax(R).$ It is easy to see that this extension is unique in this case. Also, if $R$ is $\ast$-extendible, then $\Qrmax(R)$ is its maximal left ring of quotients (the isomorphism of $R$ to the opposite ring of $R$ extends to an isomorphism of $\Qrmax(R)$ and its opposite ring). Thus, we can write $\Qmax$ for $\Qrmax(R)=\Qmax^l(R).$ If $R$ is right nonsingular (thus left nonsingular as well since it is a $\ast$-ring), then $\Qmax$ is regular and left and right self-injective. Thus, $\Qmax$ is unit-regular (see \cite[Theorem 9.29, p. 105]{Go}).

If $\ast$ is proper in a nonsingular and $\ast$-extendible ring $R$, then the extension of $\ast$ is proper in $\Qmax.$ To see that let $x^{\ast}x=0$ for some $x\in\Qmax.$ This implies that $(xr)^{\ast}xr=r^{\ast}x^{\ast}xr=0$ and so $xr=0$ for all $r$ such that $xr\in R.$ Since $I=\{r\in R|xr\in R\}$ is a dense right ideal of $R,$ $xI=0$ implies that $x=0.$ The involution in every Rickart $\ast$-ring is proper and every Rickart $\ast$-ring is nonsingular. Thus, if $R$ is $\ast$-extendible Rickart $\ast$-ring, then $\Qmax$ is $\ast$-regular.

In \cite{Handelman}, Handelman proved that a $\ast$-extendible Baer $\ast$-ring has the same projections as its maximal ring of quotients $\Qmax.$
As is the case with many other properties of Baer $\ast$-rings, sufficient conditions for $\ast$-extendibility have been given in terms of certain Baer $\ast$-ring axioms. Let $R$ be a Baer $\ast$-ring and let us look at the following axioms.
\begin{itemize}
\item[(A1)] The ring $R$ is {\em finite} if $x^{\ast}x=1$ implies
$xx^{\ast}=1$ for all $x\in R.$

\item[(A2)] The ring $R$ satisfies {\em existence of projections
(EP)-axiom:} for every $0\neq x\in R,$ there exists a self-adjoint
$y\in\{x^{\ast}x\}''$ such that $(x^{\ast}x)y^2$ is a nonzero projection;

The ring $R$ satisfies the {\em unique positive square root (UPSR)-axiom:}
for every $x\in R$ such that $x=x_1^{\ast}x_1+x_2^{\ast}x_2+\ldots
+x_n^{\ast}x_n$ for some $n$ and some $x_1,x_2,\ldots,x_n\in R$ (such
$x$ is called {\em positive}), there is a unique $y\in\{x^{\ast}x\}''$
such that $y^2=x$ and $y$ positive. Such $y$ is denoted by
$x^{1/2}.$

\item[(GC)] The ring $R$ satisfies the {\em generalized comparability} if for every two projections $p$ and $q,$ there is a central projection $c$ such that $ cp\preceq cq\;\;\mbox{ and }\;\;(1-c)q\preceq(1-c)p.$ Here $p\preceq q$ means that $p\sim r\leq q$ for some projection $r$.

\item[(LP$\sim$RP)] For every $x\in R,$ $LP(x)\sim RP(x).$
\end{itemize}

It is easy to see that (A1) is equivalent to the condition that 1 is a finite projection (a projection is said to be {\em finite} if it is not equivalent to its proper subprojection).

By \cite[Proposition 2.10]{Handelman}, a Baer $\ast$-ring with the following properties is $\ast$-extendible: (1) $R$ is finite; (2) (LP$\sim$ RP) holds; and (3) $R$ has sufficiently many projections (for every $x$ there is $y$ with $xy$ projection or, equivalently, every nonzero right (or left) ideal contains a nonzero projection). These conditions are met for finite Baer $\ast$-rings satisfying (A2) axiom: (LP$\sim$RP) follows from (A2) by \cite[Corollary on p. 131]{Be} and the existence of sufficiently many projections is guaranteed by (EP) axiom. Thus, every Baer $\ast$-ring satisfying (A1) and (A2) is $\ast$-extendible.

Let us consider the following additional axioms.

\begin{itemize}
\item[(A3)] Partial isometries are addable.

\item[(A4)] The ring $R$ is {\em symmetric}: for all $x\in R,$ $1+x^{\ast}x$ is
invertible.

\item[(A5)] There is a central element $i\in R$ such that $i^2=-1$
and $i^{\ast}=-i.$

\item[(A6)] The ring $R$ satisfies the {\em unitary spectral (US)-axiom:}
for a unitary $u\in R$ with RP$(1-u)=1,$ there is an
increasingly directed sequence of projections $p_n\in\{u\}''$ with $\sup p_n=1$
such that $(1-u)p_n$ is invertible in $p_n Rp_n$ for
every $n.$

\item[(A7)] The ring $R$ satisfies the {\em positive sum (PS)-axiom;} if
$p_n$ is orthogonal sequence of projections with supremum 1 and
$a_n\in R$ such that $a_n$ is positive (as defined in axiom (A2)) $a_n\leq p_n,$ then there is $a\in R$
such that $a p_n=a_n$ for all $n.$
\end{itemize}

Every finite $AW^{\ast}$-algebra (so a finite von Neumann algebra as well) satisfies axioms (A1) -- (A7) (\cite[Remark 1, p. 249]{Be}).

In \cite{Be}, Berberian uses (A1)--(A6) to embed a ring $R$ in a regular ring $Q$: he shows that  a Baer $\ast$-ring $R$ satisfying (A1)-- (A6) can be embedded in a regular Baer $\ast$-ring $Q$ satisfying (A1) -- (A6) such
that $R$ is $\ast$-isomorphic to a $\ast$-subring of $Q$, all projections,
unitaries and partial isometries of $Q$ are in $R$ (see \cite[Chapter 8]{Be}) and that $Q$ is
unique up to $\ast$-isomorphism (see \cite[Proposition 3, p. 235]{Be}).
The ring $Q$ is $\ast$-regular as well since it is regular and Baer (thus Rickart as well) $\ast$-ring. It is called the {\em regular ring of Baer $\ast$-ring $R$.}

The regular ring $Q$ is also the maximal and the classical ring of quotients of $R$ (\cite[Proposition 3]{Lia2}). Thus, it is left and right self-injective. Every self-injective (both left and right) and regular ring is unit-regular (see \cite[Theorem 9.29, p. 105]{Go}) so $Q$ is unit-regular as well.

If $R$ satisfies (A1) -- (A7), the ring $M_n(R)$  of $n\times n$ matrices over $R$ is a Rickart $\ast$-ring for every $n$ (by \cite[Theorem 1, p. 251]{Be}). Moreover, $M_n(R)$ is semihereditary (\cite[Corollary 5]{Lia2}). In \cite{Lia3} it is shown that $M_n(R)$ is a finite Baer $\ast$-ring for every $n$ (\cite[Theorem 4]{Lia3}). This gave an affirmative answer to the question of Berberian (\cite[Exercise 4D, p.253]{Be}) whether $M_n(R)$ is a Baer $\ast$-ring if $R$ satisfies (A1)--(A7). In \cite{Be}, Berberian uses additional two axioms (called (A8) and (A9) in \cite{Lia2}) to prove that $M_n(R)$ is finite. In \cite{Lia3}, it is shown that (A9) follows from (A1)--(A7).

Mimicking the type decomposition for von Neumann algebras, the three types of Baer $\ast$-rings have been considered. A Baer $\ast$-ring is said to be of:
\begin{itemize}
\item[(i)] {\em type $I$} if it has faithful abelian projection (a projection $p$ is {\em abelian} iff $pRp$ is abelian ring and $p$ is {\em faithful} if there are no nontrivial central idempotents $e$ with $ep = 0$);

\item[(ii)] {\em type $II$} if it has faithful finite projection but no nontrivial abelian projections; and

\item[(iii)] {\em type $III$} if there is no nontrivial finite projections.
\end{itemize}
A Baer $\ast$-ring has a unique decomposition into three Baer $\ast$-rings of the three different types.

Moreover, a finer type decomposition of types $I$ and $II$ is possible. A Baer $\ast$-ring is said to be of:
\begin{itemize}
\item[-] {\em type $I_{f}$} if it is of type $I$ and finite;

\item[-] {\em type $I_{\infty}$} if it is of type $I$ and  0 is the only abelian projection;

\item[-] {\em type $II_1$} if it is finite and 0 is the only finite central projection;

\item[-] {\em type $II_{\infty}$} if 0 is the only abelian and the only finite central projection.
\end{itemize}
Every Baer $\ast$-ring has a unique decomposition into five Baer $\ast$-rings of the following types: $I_f,$ $I_{\infty},$ $II_1,$ $II_{\infty}$ and $III$. For details on type decompositions see \cite[Theorems 2 and 3, p. 94]{Be}.

Type $I_{f}$ rings can be further classified. Let $R$ be of type $I_{f}$ such that there is a positive integer  $n$ and $n$ equivalent orthogonal abelian projections that add to 1 (called homogeneous partition of 1 of order $n$). In this case $R$ is said to be of type $I_n.$ If $R$ is of type $I_{f}$ and has (GC), it can be uniquely decomposed into a $C^*$-sum of types $I_n$ (\cite[Theorem 2, p. 115]{Be}). If $R$ has (GC) and is of type $I_n$ for some $n,$ then $R$ is $\ast$-isomorphic to $M_n(A)$ where $A$ is an abelian Baer $\ast$-ring uniquely determined up to $\ast$-isomorphism (\cite[Proposition 2, p. 112]{Be}). This gives us the following result.

\begin{thm} If a Baer $\ast$-ring $R$ is of type $I_n$ and satisfies (A2), it is almost $\ast$-clean. If, in addition, $R$ is regular, then it is $\ast$-clean.
\label{type_In}
\end{thm}
\begin{pf}
Note that $R$ satisfies (GC) since (A2) holds (\cite[Theorem 1, p. 80]{Be}). Thus, $R$ is $\ast$-isomorphic to $M_n(A)$ where $A$ is an abelian Baer $\ast$-ring uniquely determined up to $\ast$-isomorphism (\cite[Proposition 2, p. 112]{Be}). Moreover, if $p_1, p_2,\ldots, p_n$ are orthogonal abelian projections forming a homogeneous partition of 1 of order $n,$ then $A$ is $\ast$-isomorphic to the corner $p_1 Rp_1$ (and $p_i Rp_i$ for any $i=1,\ldots, n$).

A corner $pRp$ of a $\ast$-ring $R$ where $p$ is a projection, preserves the following properties: Rickart  $\ast$-ring, Baer  $\ast$-ring (see \cite[Theorem 4, p. 6, and p. 30]{Kaplansky}), regular (if $x\in pRp$ then $x=xyx$ implies $x=xyx=(xp)y(px)=x(pyp)x$), $\ast$-regular (if $R$ has a proper involution then it is easy to see the involution in $pRp$ inherited from $R$ is also proper), and regular and right and left self-injective (see \cite[Proposition 13.7, p. 98]{Go}). If $R$ is a Baer $\ast$-ring with $Q=\Qmax(R)$ and $p$ is a projection in $R$, then $pQp$ is the maximal ring of quotients of $pRp.$ This is because $R$ is semiprime (every Baer $\ast$-ring is semiprime) and nonsingular (every Rickart $\ast$-ring is nonsingular) so $pQp$ is an essential extension of $pRp$ that is right self-injective and so it is equal to $\Qmax(pRp)$ (see \cite[Proposition 13.39]{Lam}). This last result is also in \cite[Proposition 0.2, p. 135]{Cail_Ren}.

Thus, the maximal ring of quotients of the abelian Baer $\ast$-ring $A\cong p_1Rp_1$ is $p_1Qp_1.$ The ring $p_1Qp_1$ is an abelian $\ast$-regular ring so it is $\ast$-clean by  Proposition \ref{clean_regular}. The partition  $p_1, p_2,\ldots, p_n$  is a homogeneous partition of order $n$ for $Q$ as well. Since $p_iQp_i\cong p_1Qp_1\cong\Qmax(A)$ is $\ast$-clean for every $i,$ $Q$ is $\ast$-clean by Proposition \ref{*-clean_matrices}.

The ring $R$ satisfies (A1) and (A2) and so it is $\ast$-extendible. Thus, the projections of $R$ and $Q$ are the same. Then $R$ is almost $\ast$-clean by Proposition \ref{embedding_proposition}.
\end{pf}

Now we can prove the almost $\ast$-cleanness of rings of type $I_f.$ First we prove the following proposition.

\begin{prop} Let $R$ be a Baer $\ast$-ring that satisfies (A1) -- (A6). If there are central orthogonal projections $p_n$ with supremum 1 and $p_nR$ are either of type $I_n$ or equal to 0, then $R$ is almost $\ast$-clean. If $R$ is also regular, then (A6) is not needed, $R$ is isomorphic to a direct product of rings $p_nR$ and $R$ is $\ast$-clean.
\label{Baer_C*-sums}
\end{prop}
\begin{pf}
As we noted before, axioms (A1)--(A6) are sufficient to guarantee that $R$ embeds in a $\ast$-regular ring $Q$ with same projections as $R$ satisfying (A1)--(A6) (for more details see \cite[Section 52]{Be}). It is easy to see that $p_nR$ is also a Baer $\ast$-ring that satisfies (A1)--(A6) for every $n$ (here it is essential that the projections $p_n$ are central). By uniqueness of regular ring of Baer $\ast$-ring satisfying (A1)--(A6) (see \cite[Proposition 3, p. 235]{Be}), it is easy to see that $p_nQ$ is the regular ring of $p_nR$ for all $n$ and thus that $p_nR$ and $p_nQ$ have the same projections for every $n$. Since $p_nR$ and $p_nQ$ have the same projections and $p_nR$ is either trivial or of type $I_n,$ $p_nQ$ is also either trivial or of type $I_n$ as well. Then $p_nQ$ is $\ast$-clean by Theorem \ref{type_In} for every $n$.

The crucial ingredient in the rest of the proof is the following result given in \cite[Theorem 2, p. 237]{Be}.
\begin{itemize}
\item[(R)] If $R$ is a Baer $\ast$-ring that satisfies (A1) -- (A6) and there are central orthogonal projections $p_n$ with supremum 1, then the regular ring $Q$ of $R$ is $\ast$-isomorphic to the direct product of rings $p_nQ$ via $x\mapsto (p_nx).$
\end{itemize}

Since a direct product of $\ast$-clean rings is $\ast$-clean (easy to see), the regular ring $Q$ of ring $R$ is $\ast$-clean by (R). Then $R$ is almost $\ast$-clean by Proposition \ref{embedding_proposition}.

To prove the last sentence of the statement of this proposition, note that if $R$ is regular and satisfies (A1)--(A5), then (A6) holds as well (see \cite[Exercise 4A, p. 247]{Be}) and $R$ is equal to its regular ring. Then $R$ is isomorphic to a direct product of rings $p_nR$ by (R). In this case, the regular rings $p_nR$ are $\ast$-clean by Theorem \ref{type_In} and so $R,$ a direct product of $\ast$-clean rings, is $\ast$-clean as well.
\end{pf}

\begin{thm}
If a Baer $\ast$-ring $R$ is of type $I_{f}$ and satisfies (A2) -- (A6), it is almost $\ast$-clean. If $R$ is also regular, it is $\ast$-clean.
\label{type_If}
\end{thm}
\begin{pf}
Every Baer $\ast$-ring of type $I_f$ that satisfies (A2) has a sequence of orthogonal central projections $p_n$ with supremum 1 such that $p_nR$ is either 0 or of type $I_n$ (see \cite[Theorem 2, p. 115]{Be}). Condition (PC) assumed in \cite[Theorem 2, p. 115]{Be} follows from (GC) (by \cite[Proposition 2, p. 78]{Be}) and (GC) follows from (A2) (\cite[Theorem 1, p. 80]{Be}). Thus every type $I_f$ Baer $\ast$-ring that satisfies (A2)--(A6), satisfies the assumptions of Proposition \ref{Baer_C*-sums}. Thus $R$ is almost $\ast$-clean. If, in addition, $R$ is regular, $R$ is $\ast$-clean also by Proposition \ref{Baer_C*-sums}.
\end{pf}

Let us conclude this section with the observation that Baer $\ast$-rings with (A1)--(A6) are morphic exactly when they are regular (and this happens exactly when the rings are their own regular rings). Recall that a ring is right morphic if $\ann_r(x)\cong R/xR$ for every $x\in R.$ The following five conditions are equivalent (see  \cite[Exercise 19A, p. 270]{Lam} and \cite[Corollary 3.16]{Zhu_Ding}): (1) $R$ is unit-regular, (2) $R$ is regular and right morphic, (3) $R$ is regular and left morphic, (4) $R$ is right Rickart and left morphic, and (5) $R$ is left Rickart and right morphic. Thus, if $R$ is both left and right Rickart, the conditions of being unit-regular, left morphic and right morphic are equivalent.

In \cite{Nicholson_Sanchez},  it is shown that every right morphic ring is left principally injective (also called left P-injective) meaning that every map $Rx\rightarrow R$ extends to $R$ (\cite[Theorem 24]{Nicholson_Sanchez}). If $R$ is both left and right Rickart, the conditions of being regular, left P-injective and right P-injective are equivalent. This follows from \cite[Proposition 27]{Nicholson_Sanchez} and \cite[Corollary 3.15]{Zhu_Ding}.

In \cite{Camillo_Nicholson}, a weaker notion of right morphic ring, called right quasi-morphic, is considered: for every $x\in R,$ there is $y$ and $z$ such that $xR = \ann_r(y)$ and $\ann_r(x) = zR$. It is shown that the results of P-injectivity from \cite{Nicholson_Sanchez} hold if ``morphic'' is replaced by ``quasi-morphic''.

\begin{cor} Let $R$ be a Baer $\ast$-ring that satisfies  (A1)--(A6). If  $Q$ is the regular ring of $R$, then the following conditions are equivalent:

\begin{tabular}{ll}
(1) $R=Q.$ & \\
(2) $R$ is regular. &  (3) $R$ is unit-regular. \\
(4) $R$ is left morphic. & (5)  $R$ is right morphic.\\
(6) $R$ is left quasi-morphic. & (7)  $R$ is right quasi-morphic.\\
\end{tabular}
\label{morphic}
\end{cor}
\begin{pf}
(1) $\Leftrightarrow$ (2).  If $R=Q$ then $R$ is regular since $Q$ is. Conversely, if $R$ is regular, then $R=Q$ by uniqueness of the regular ring (\cite[Proposition 3, p. 235]{Be}).

(2) $\Leftrightarrow$ (3).  If $R$ is regular, then $R=Q$ implies that $R$ is its own maximal (left and right) ring of quotients (see  \cite[Proposition 3]{Lia2} and note that the proof uses just (A1)--(A6) and not (A7)). Thus, $R$ is regular and (left and right) self-injective. Then $R$ is unit-regular (see \cite[Theorem 9.29, p. 105]{Go}). The converse (3) $\Rightarrow$ (2) always holds.

(3) $\Leftrightarrow$ (4) $\Leftrightarrow$ (5). As observed above using results from \cite{Zhu_Ding}, this is true for every Rickart $\ast$-rings.

(6) $\Rightarrow$ (2). As observed above using result from \cite{Zhu_Ding}, if $R$ is both left and right Rickart, the conditions of being regular, left P-injective and right P-injective are equivalent. Since left quasi-morphic implies right P-injective by \cite[Lemma 3]{Nicholson_Sanchez}, left quasi-morphic implies regular.

By symmetry, (7) $\Rightarrow$ (2) as well. The condition (3) implies (4) and (5) which imply (6) and (7) respectively. Hence, all seven conditions are equivalent.
\end{pf}

Note that in general the following hold: (2) is weaker than (3), (4) and (5) are weaker than (3), and neither (4) nor (5) implies the other. The examples of rings illustrating these claims can be found in \cite{Nicholson_Sanchez}. Also, (6) and (7) are weaker than (3) -- (5) and  neither (6) nor (7) implies the other. The examples of rings illustrating these properties can be found in \cite{Camillo_Nicholson}.

\section{$\ast$-clean and almost $\ast$-clean von Neumann algebras}
\label{section_VNAs}

In this section, we turn to $AW^{\ast}$ and von Neumann algebras. First, let us note that all $AW^{\ast}$-algebras satisfy (A2)--(A7) (see \cite[pp. 48, 70, 129, 233 and 244]{Be} and \cite[p. 327]{KadRing}). In particular, all finite $AW^{\ast}$-algebras satisfy (A1)--(A7). Thus we have the following.
\begin{cor} An $AW^{\ast}$-algebra of type $I_{f}$ (in particular a von Neumann algebra of type $I_f$) is almost $\ast$-clean. If it is also regular, it is $\ast$-clean.
\label{VNA_clean}
\end{cor}

It is interesting to note that one of the key arguments in the proof of Proposition \ref{Baer_C*-sums} is that a direct product of $\ast$-clean rings is $\ast$-clean. For $C^{\ast}$-algebras, however, the concepts of direct product is not suitable concept for ``product'' since the coordinate-wise norm might fail to be complete. A direct sum of algebras, on the other hand, might fail to be a unitary algebra. Because of this, a more appropriate concept of a product is considered.

If $\{A_i\}$ is a family of $C^{\ast}$-algebras, its {\em $C^{\ast}$-sum} is defined to be the $\ast$-subalgebra of the direct product $\prod A_i$ containing all elements $(a_i)$ such that $\sup ||a_i||$ is finite. Then the supremum norm is complete and this is a $C^{\ast}$-algebra (see \cite[Section 10]{Be}). In \cite{KadRing}, this algebra is denoted by $\sum\oplus\, A_i.$ The properties of being Rickart and Baer are preserved by $C^{\ast}$-sums of $C^{\ast}$-algebras and so a $C^{\ast}$-sum of Rickart $C^{\ast}$ and $AW^{\ast}$-algebras remain Rickart $C^{\ast}$ and $AW^{\ast}$ respectively (\cite[Proposition 1, p. 52]{Be}).

When working with families of central orthogonal projections with supremum 1, $C^{\ast}$-sums provide an appropriate framework  because of the following result (\cite[Proposition 2, p. 53]{Be}).
\begin{itemize}
\item[($\sum\oplus$)]  If $\{p_i\}$ is an orthogonal family of central projections with supremum 1 in an   $AW^{\ast}$-algebra $A,$ then $A$ is $\ast$-isomorphic to $\sum\oplus\, p_iA$ via $x\mapsto (p_ix).$
\end{itemize}

\begin{prop} If $A_i$ are $I_{f}$ type $AW^{\ast}$-algebras, then $\sum\oplus\, A_i$ is almost $\ast$-clean. If $A_i$ are also regular,  $\sum\oplus\, A_i$ is $\ast$-clean.
\label{AW_C*-sums}
\end{prop}
\begin{pf}
The proof follows the idea of the proof of Proposition \ref{Baer_C*-sums}.
Let $A$ denote $\sum\oplus\, A_i.$ It is easy to see that $A$ is a finite $AW^{\ast}$-algebra as well so it has the regular ring. Let $Q_i$ be the regular ring of $A_i$ and $Q$ the regular ring of $A.$ Since the projections in $Q$ and $A$ are the same, there are central orthogonal projections $p_i$ in $Q$ with supremum 1 such that $p_iQ\cong Q_i$ by uniqueness of the regular ring. So $Q$ is $\ast$-isomorphic to the $C^{\ast}$-sum of $Q_i.$ On the other hand, since $Q$ is regular, it is $\ast$-isomorphic to direct product of $Q_i$ by result (R) quoted in the proof of Proposition \ref{Baer_C*-sums}. The rings $Q_i$ are $\ast$-clean by Corollary \ref{VNA_clean} (rings $Q_i$ are of type $I_f$ because rings $A_i$ are of type $I_f$ and the projections in $Q_i$ and $A_i$ are the same). Thus, $Q$ is $\ast$-clean since it is a direct product of $\ast$-clean rings. Since the projections in $A$ and $Q$ are the same, $A$ is almost $\ast$-clean by Proposition \ref{embedding_proposition}.

If $A_i$ are regular, $A_i=Q_i$ and so $A_i$ are $\ast$-clean. The $C^{\ast}$-sum $\sum\oplus\, A_i$ is $\ast$-isomorphic to the direct product so it is $\ast$-clean as well.
\end{pf}

We shall exhibit an example of a von Neumann algebra that is $\ast$-clean and not of type $I_f$. An open question which still remains is whether there are such examples of type $II_1$ (and the remaining three types as well). The consideration of factors of type $II_1$ would provide an important step in proving (or disproving) the (almost) cleanness of algebras of type $II_1$. To motivate further investigation in this direction, we finish the paper with some examples and some questions.

Let us consider the class of group von Neumann algebras. This class provides us with some concrete examples of finite von Neumann algebras.
Let $G$ be a group. The complex group ring $\Cset G$ is a
pre-Hilbert space with an inner product and an involution given by
\[ \langle \;\sum_{g\in G}
a_g g, \sum_{h\in G} b_h h\;\rangle = \sum_{g\in
G}a_g\overline{b_g}\;\;\mbox{ and }\;\;\left(\sum_{g\in G} a_g g\right)^{\ast} = \sum_{g\in G} \overline{a_g}
g^{-1}.\]

The Hilbert space completion of $\Cset G$ is $\lge,$ the space of square summable complex valued functions over the group $G$: \[ \lge = \{\;\sum_{g\in G} a_g g\;|\;\sum_{g\in G}
|a_g|^2 < \infty\;\}.\] The involution from $\Cset G$ extends to $\lge.$

{\em The group von Neumann algebra} $\vng$ is the space of
$G$-equivariant bounded operators from $\lge$ to itself: \[\vng =
\{\; f\in {\B}(\lge)\;|\;f(xg)=f(x)g\mbox{ for all }g\in G\mbox{
and }x\in\lge\;\}.\] The algebra $\vng$ is a von Neumann algebra on $H = \lge$
since it is the weak closure of $\Cset G$ in ${\B}(\lge)$ so it is
a $\ast$-subalgebra of ${\B}(\lge)$ which is weakly closed
(see \cite[Example 9.7]{Lu5} for details). The algebra $\vng$ is finite since it has a normal, faithful trace
tr$_{\A}(f) = \langle f(1), 1 \rangle$ (every finite von Neumann algebra has such a trace and, conversely, a von Neumann algebra with such a trace is finite, for details see \cite[Section 9]{Lu5}).
The regular ring of $\vng$ is the algebra of affiliated operators $\ug$ (see \cite[Section 8]{Lu5}).

Since $\vng$ is finite, just types $I_{f}$ and $II_1$ are possible. Recall that a group is said to be virtually abelian if it has an abelian subgroup of finite index. The types of $\vng$ are classified according to the properties of group $G$ as follows:
\begin{itemize}
\item[(i)] The algebra $\vng$ is of type $I_{f}$ iff $G$ is virtually abelian. If $G$ is finitely generated, $\vng$ is of type $II_1$ iff $G$ is not virtually abelian.

\item[(ii)] The algebra $\vng$ is of type $II_1$ iff $G_f$ has infinite index ($G_f$ is the normal subgroup of $G$ of elements with finitely many elements in the conjugacy class or, equivalently, of elements whose centralizer has finite index).

\item[(iii)] The algebra $\vng$ is a factor (its center is isomorphic to $\Cset \{1\}$) iff $G_f$ is trivial.
\end{itemize}
\cite[Lemma 9.4, p. 337]{Lu5} contains more details. Moreover, $\vng$ is semisimple iff $G$ is finite (see \cite[Exercise 9.11, p. 367]{Lu5}). Note also that an $AW^{\ast}$-algebra  is commutative iff it is abelian (\cite[Example 2, p. 90]{Be}). Thus, $G$ is an abelian group iff $\vng$ is an abelian $\ast$-ring. Now let us consider the following examples.

\begin{exmp} {\em
\begin{enumerate}
\item Finite abelian groups give rise to commutative and semi\-simple group von Neumann algebras. Finite and not abelian groups give rise to non-commutative (nor reduced, nor Armendariz) semisimple group von Neumann algebras.

\item Let $G=\Zset.$ Then $\vng$ is commutative and not semisimple (since $G$ is infinite). Moreover, $\vng$ is not regular: $\vng$ can be identified with the space of (equivalence classes of) essentially bounded measurable complex function on the unit circle (see \cite[Example 1.4]{Lu5}). Its regular ring, the algebra of affiliated operators $\ug$ of $\vng$, can be identified with the space of (equivalence classes of) {\em all} measurable complex functions on a unit circle (see \cite[Example 8.11]{Lu5}). Thus $\vng\neq\ug$ and so $\vng$ is not regular by Corollary \ref{morphic}.

By Corollary \ref{morphic}, $\vng$ is not morphic (nor quasi-morphic) since it is not regular. On the other hand, $\ug$ is strongly $\ast$-clean by Proposition \ref{clean_regular} and $\vng$ is almost strongly $\ast$-clean by Proposition \ref{embedding_proposition}. Thus, this is an almost strongly clean Baer $\ast$-ring that is not morphic and not regular.

Note also that the group ring  $\Cset G$ is not clean by \cite[Proposition 2.7]{Warren_RG} (stating that if $R$ and $G$ are commutative and the group ring $RG$ is a clean ring, then $R$ is clean and $G$ is a torsion group.)
Since $\vng$ is almost clean, we have the following situation: a non-clean ring is dense (in the topological sense, since $\vng$ is the closure of $\Cset G$ in the weak topology) in an almost clean ring which is dense in a clean ring ($\ug$).

\item Let $G=\Zset\oplus D_3$ where $D_3$ is the dihedral group (group of symmetries) of an equilateral triangle.
This is an example of non-regular (so non-morphic as well) and non-abelian ring. The group $G$ is virtually abelian and so $\vng$ is of type $I_f$. Thus $\vng$ is almost $\ast$-clean by Corollary \ref{VNA_clean}.
\end{enumerate}

Note that the algebras in all of the above examples are $I_f$ type (because all the groups considered are virtually abelian). Now let us consider some non-virtually abelian groups.

Let  $\{G_{\alpha}\}_{\alpha\in \varLambda}$ be a family of groups. Then
\[\N(\prod G_{\alpha})\cong \sum\oplus\, \N(G_{\alpha})\leq \prod \N(G_{\alpha})\]
The last inclusion could be strict if the family of groups is infinite and the algebras $\N(G_{\alpha})$ are not regular (e.g. every $G_{\alpha}$ is $\Zset$). This is because $\prod \N(G_{\alpha})$ does not have to be complete and  $\sum\oplus\, \N(G_{\alpha})$ is always complete. The first two algebras are isomorphic since the elements $p_{\alpha}=(\delta_{\alpha\beta}1_{\beta})_{\beta\in \varLambda}$ where $1_{\alpha}$ is the identity element in $G_{\alpha},$ are central orthogonal projection with supremum 1 in $\N(\prod G_{\alpha}).$ Thus, $\N(\prod G_{\alpha})$ is $\ast$-isomorphic to $\sum\oplus\, p_{\alpha}\N(\prod G_{\alpha})=\sum\oplus\, \N(G_{\alpha})$ by ($\sum\oplus$) (\cite[Proposition 2, p. 53]{Be}).

If all groups $G_{\alpha}$ are finite, the algebras $\N(G_{\alpha})$ are regular and so $\sum\oplus \N(G_{\alpha})=$ $\prod \N(G_{\alpha})$ by Proposition \ref{Baer_C*-sums} and by  ($\sum\oplus$). Note that we have equality and not just isomorphisms because the isomorphism from Proposition \ref{Baer_C*-sums} and from  ($\sum\oplus$) that maps $x=(x_{\alpha})$ onto $(p_{\alpha}x)=(x_{\alpha}),$ is identity in this case. Thus, $\sum\oplus\, \N(G_{\alpha})=\prod \N(G_{\alpha})\cong\N(\prod G_{\alpha})$ and $\N(\prod G_{\alpha})$ is $\ast$-clean.

\begin{itemize}
\item[(4)]  Let $\{G_n\}_{n=1}^{\infty}$ be a countably infinite family of finite groups and let $G=\prod G_n.$ If $G_n$ are abelian, $\vng$ is an example of a type $I_f$ ($I_1$ in fact) regular ring that is not semisimple. If infinitely many groups $G_n$ are not abelian (for example, we can take $G_n=D_3$ for every $n$), then $\vng$ is a regular and $\ast$-clean (by Proposition \ref{AW_C*-sums}) ring that is not of type $I_f$ (since $G$ is not virtually abelian).
\end{itemize}

This last example shows that there are clean Baer $\ast$-rings outside of type $I_f.$

Finally, let us consider some groups of type $II_1.$

\begin{itemize}
\item[(5)] Let $G=\Zset\ast\Zset$ be the free group on two generators.  Then $G_f$ is trivial and has infinite index and so $\vng$ is a factor of type $II_1.$ Is it (almost) clean?

Note also that we can easily construct non-$\ast$-isomorphic factors of type $II_1$. For example, let $\varPi$ be the group of permutations on a countably infinite set that leave fixed all but finitely many elements (\cite[Example 6.7.7]{KadRing2}). Then $\N(\Zset\ast\Zset)$ and $\N(\varPi)$ are both factors of type $II_1$ but are not $\ast$-isomorphic (see \cite[Theorem 6.7.8]{KadRing2}).
\end{itemize} }
\label{Group_Example}
\end{exmp}

\section{Questions}
\label{section_questions}

We conclude the paper with a list of some questions.
\begin{enumerate}
\item Find an example of a $\ast$-ring that is clean but not $\ast$-clean. By Lemma \ref{abelian_lemma}, such example cannot be found in the class of abelian Rickart $\ast$-rings.

\item Which type $I_f$ $AW^{\ast}$-algebras (and Baer $\ast$-rings with (A1)--(A6)) are clean? We know that regular type $I_f$ $AW^{\ast}$-algebras are clean. Are there examples outside of this class?

\item Are $AW^{\ast}$-algebras (and Baer $\ast$-rings with (A1)--(A6)) of type $II_1$ (almost) clean? A possible first step in this consideration might be to determine if the group von Neumann algebra of the free group on two generators is almost clean. Then, consideration of other types could follow leading to a complete answer to Lam's question on cleanness of von Neumann algebras.

\item Do Theorem \ref{type_In}, Theorem \ref{type_If} and Corollary \ref{VNA_clean} hold if ``almost clean'' is replaced by ``almost strongly clean''? In case of group von Neumann algebras, are all type $I_f$ group von Neumann algebras almost strongly clean (i.e. if $G$ is virtually abelian, is $\vng$ almost strongly clean)?

Also, it would be interesting to see which results that hold for strongly clean rings hold for strongly $\ast$-clean rings as well. In addition, in \cite{Nich2} Nicholson asked if a unit-regular ring is strongly clean. We ask if a unit-regular and $\ast$-regular ring is strongly $\ast$-clean.
\end{enumerate}

\end{document}